\documentclass[12pt]{article}
\usepackage{ams}
\newtheorem{theorem}{Theorem}
\newtheorem{definition}{Definition}
\newtheorem{lemma}{Lemma}
\newtheorem{suggestion}{Suggestion}

\begin{document}

\title{The  Diophantine  equations \\
$ x^{n}_{1}  +x^{n}_{2}  +\ldots+x^{n}_{r_{1}}= y ^{n}_{1}
+y^{n}_{2}  +\ldots+y^{n}_{r_{2}}   $}

\author
{Michael A. Ivanov \\
Physics Dept.,\\
Belarus State University of Informatics and Radioelectronics, \\
6 P. Brovka Street,  BY 220027, Minsk, Republic of Belarus.\\
E-mail: ivanovma@gw.bsuir.unibel.by.}
\date{}
\maketitle

\begin{abstract}

               The aim of this paper is to prove the possibility of
            linearization of such equations by means of introduction of
            new variables. For $n=2$ such a procedure is well known, when
            new variables are components of spinors and they are widely
            used in mathematical physics. For example, parametrization of
            Pythagoras threes $a^{2} +b^{2}$ , $a^{2} -b^{2}$ , $2ab$ may
            be cited as an  example in number theory where two independent
            variables  form  a spinor which can be obtained by solution of
            a system of two linear equations. \par
               We also investigate the combinatorial estimate  for  the
            smallest sum $r(n)=r _{1}+r_{2} -1 $ for solvable equations of such a
            type as $r(n) \leq 2n+1$ (recently the better one with $r(n) \leq2n-1$ was
            received by L. Habsieger (J. of Number Theory 45 (1993) 92)).
            Apart from that we consider two conjectures about $r(n)$ and
            particular solutions for $n \leq11$ which were found with the help
            of the algorithm that is not connected with linearization.
\end{abstract}

\section[ 1.]{ Introduction                     }
               A great interest was displayed to the particular case
            of such equations - $x ^{n}+y^{n} =z ^{n}$  \cite{1,2}, and it seems rather
            strange, that the general case of equality of sums of the
            same powers
\begin{equation}
                          \sum_{1} ^{r_{1}} x_{i}^{n}  =\sum_{1}^{r_{2}}    y_{j}^{n}
\end{equation}
            was left in the shadow, although for large  $r_{1} +r_{2}$   such
            equations are solvable  and  for  small  $r _{1}+r_{2}$   a  set  of
            assertions of the  type of Fermat's last  theorem  can  be
            formulated for them: $r(n)>R(n) \mid n>n_{0} $, where $R(n)$ and $n _{0}$
            are fixed, $r(n)$ is the smallest sum $r _{1}+r_{2} -1$ for which Eqs.(1)
            have at least a non-trivial solution.
            \par   A geometrical approach \cite{3} was used by G. Faltings to prove
            Mordell's conjecture \cite{4}, that gave a new push to the
            investigations of Fermat's equations. On the other hand,
            an interesting algebraic fact, concerning (1), was recently
            discovered : they can be linearized by introducing new
            variables \cite{5,6,7,8} with the help of which assumed solutions of
            Eqs.(1) can be parametrized. A different approach to
            linearization is described in \cite{9}.
            \par   In this paper we give a constructive proof of the
            theorem about the possibility of linearization of (1) (when
            n is a prime integer), the technique of which consists in
            adding a system of linear equations in which some new variables
            enter (see \cite{7}):

\begin{theorem}
            A. If we take the equation
\begin{equation}
                   \sum_{1}^{k}  x_{i}^{n}  = y^{n} ,  k \geq 2,
\end{equation}
           a linear matrix equation
\begin{equation}
                          \sum_{1}^{k}   x_{i} A_{i}\Psi  = yE\Psi,
\end{equation}
           can be juxtaposed with it, where $A_{i}$  are the square matrices
           of order  $mn$, for which  the following condition is necessary
\begin{equation}
     \forall \lbrace x_{i} \rbrace:  ~~ {( \sum_{1}^{k} x_{i} A _{i})}^{n}  =
     E \sum_{1}^{k} x_{1} ^{n}  ,
\end{equation}
           where $E$ is the unit matrix, $\Psi $   is the column-vector of new
           variables ${( \Psi_{1},\ldots, \Psi_{mn}  )}^{T} $.
           B. If (2) is fulfilled, the determinant of the system (3) with respect to unknowns  $\Psi_{i}$ is equal to zero.
\end{theorem}
\par The same matrices, which are built to linearize (2),
permit to linearize  (1), juxtaposing the equation
\begin{equation}
  \sum_{1}^{r_{1}} x_{i} A_{i} \Psi  = \sum_{1}^{r_{2}}  y_{j} B_{j}\Psi ,
\end{equation}
          with (1), where $B _{j}=\varepsilon_{j} A_{r_{1}+j},~
           \varepsilon_{j}^{n}  =-1,~ r_{1} +r _{2}=k$.
            \par
            For composite powers $n$ the  linearization  of  (1)  can  be
            carried out step by step for all prime  factors  of $ n$ using
            Theorem 1 (see also \cite{8}).
\par
The procedure of linearization is generalization of the algebraic
part of the Dirac procedure \cite{10} to the equations
$$   \varepsilon = \sum_{m=1}^{N} {(p_{m1}^{n} +p_{m2}^{n} +
\ldots+p_{mr}^{n}) }^{1/l},  $$
            and apart from that the particular case of $N=1,~l=n$ is considered in
            Theorem 1 \cite{7}. Special case of $N=l=n=2,~ r=4$, concerns
            the physical models  of  composite  particles \cite{11,12} .
            A more complicated case of it is considered in \cite{8}.
            \par   In this paper, we also consider the parametrization of
            solutions of (1) for $n=3$. Furthermore we add some remarks
            concerning solvability of Eqs.(1) and some estimates of $r(n) $
            for small $n$.
            \par The matrices $A_{i}$ , which will be built below, are the permutation
            matrices with the elements  $\zeta^{m}$, where $\zeta$ is a primitive
            root of degree $n$ of unity.

\section[ 2.]{ An example of linearization, definitions, proof of lemmas,
auxiliary and main theorems} Let us consider the equation
$x_{1}^{2}  +x_{2}^{2}  =y^{2} $ as the  simplest example and
juxtapose the linear matrix equation
           $ x_{1} A _{1}\Psi + x_{2} A _{2}\Psi = yE\Psi $ with it,
           where
$$
A_{1}=\left(  \begin{array}{cc}
                        1&0\\
                         0&-1
                         \end{array}\right),~~
A_{2}=\left(  \begin{array}{cc}
                         0&1\\
                         1&-0
                         \end{array}\right)
$$
         and $ \Psi= {( \Psi_{1},~\Psi_{2} )}^{T} $. The solutions of this  linear  equation
        have the form $x_{i} =(y/\delta)\delta_{i},~
        \delta_{1}=\Psi_{1}^{2}-\Psi_{2}^{2},~
        \delta_{2}=2 \Psi_{1}^{2}\Psi_{2}^{2},~
        \delta=\Psi_{1}^{2}+\Psi_{2}^{2}$. If we suppose that $y=\delta$,
        then one gets all Pythagoras threes for
            integers $\Psi_{i}$ .
\begin{definition}   Such a combination as
\begin{equation}
                     (ABC\ldots)_{+}\equiv  \sum_{P}   ABC \ldots,
\end{equation}
            where the summation has been made over all different permutations of
            $A,B,C,\ldots$, is the generalized anticommutator
           of matrices $A,B,C,\ldots$.
\end{definition}
               For  instance: $ {(AB)}_{+} =AB+BA$,  but $ {(A ) ^{2}}_{+}=A^{2} $,
\begin{equation}
  (A ^{k}BC\ldots)_{+} =(1/k!)(A_{1}\ldots A_{k} BC\ldots)  \vert  _{A_{i}=A}.
\end{equation}
\begin{definition} A set of $k$ matrices $A_{i}$  for some $n$ is
           a concerted one, if
$$\forall \{ n_{1} ,\ldots,n_{k}\vert~ 0\leq n_{i} <n, n_{1} +n_{2}
+\ldots+n_{k} =n \}:$$
\begin{equation}
(A _{1}^{n_{1}}  A _{2}^{n_{2}} \ldots A _{k}^{n_{k}}  )_{+} =0.
\end{equation}
\end{definition}
\begin{suggestion}  If $k$ matrices $A_{i}$  are concerted and $A _{i}^{n} =E$, the
           condition (4) is fulfilled.
\end{suggestion}

\begin{definition}  A set of positions  of  non-zero  elements is
           a structure of a permutation matrix.
\end{definition}
               To prove the part A of Theorem 1, we need two lemmas:
\begin{lemma}  There exists a concerted pair of $n\times n$ matrices.
\end{lemma}
\begin{lemma} $ m+1 $ concerted matrices can be built from $m$ ones,
           $m>1.$
\end{lemma}

               The order of matrices $A_{i}$  can be reduced thanks to Lemma 3:
\begin{lemma} There exist three concerted $n\times n$ matrices.
\end{lemma}
            (There are only three matrices of this kind \cite{8}).
             \par  Lemma 2 is proved by the author and Lemmas 1 and 3 are
            proved for the prime $n>2 $ (see \cite{8} for a more general case).
             \par {\it Proof of  lemmas.}  Let $S_{n}\equiv\{ 0,1,\ldots ,n-1 \}$ be  the  full
            system of residue classes modulo $n$. Then we can speak about
            the natural numbers of some set $S$ that they are distributed
            uniformly over the  classes of  $S_{n}$ , if the same number of
            elements of $S$ belongs to every class of  $S_{n}$. Let $S ^{m}_{nj} \subset  S_{n} $ be some
            subset of  $S_{n}$, which includes $m$ elements, with a number $j$, where
            $j=1,~C _{n}^{m}$ ($C_{n}^{m}$  is  the combination number), then
            $S _{nj}^{m}\not= S _{ni}^{m} \vert~ i\not= j$; let
            us denote $ b _{j}^{m}\equiv \sum_{1}^{m}  k_{i},~ k_{i}\in S_{nj}^{m} $.
            To prove Lemmas 1 and 3 for prime $n>2$, we shall use
\begin{theorem} For every $m, ~1\leq m<n$, the numbers $b_{j}^{m}$  are
distributed uniformly over the  classes of $S_{n}$.
\end{theorem}
              {\it Proof for the prime $n>2$ (by induction).}
            \par  1. The  theorem  is evident for $m=1.$ \par 2. Let for $\forall m',~
            m'\leq m-1$, the numbers $b_{j}^{m'}$
            be distributed uniformly over residue  classes. Then we prove,
            that $b_{j}^{m}$  are distributed uniformly. We have $\{(a+k_{i} )\vert ~
            k_{i}\in S_{n} \}= S_{n}$
            for integer $a$.  By the assumption of induction $b_{j}^{m-1}$    are
            distributed uniformly; if we add all possible $k _{l}\in S _{n}$ to each
            sum $b _{j}^{m-1}$, these  numbers will be distributed  uniformly too.
            To get $\{ b_{j}^{m}\}$, we should exclude $m-1$ numbers from $n$ ones. Then
            the excluded numbers will have the following view:
$$        a_{j} =2k_{j} +b_{l}^{m-2}\vert~    k _{j}\in S_{nl}^{m-2},$$
            with $\{2k _{j}(mod~ n)\} = S_{n}$  for the prime $n>2$, where one part of
            numbers $a_{j}$  is
$$  \{ a _{j}\vert ~ \exists  k_{i}\equiv 2k_{j} (mod ~ n),~ k_{i} \in S_{nl}^{m-2}  \}   = \{b _{j}^{m-1}\}   $$
            which is distributed  uniformly; the other part of $a_{j}$  is

$$      \{ a _{j}\vert~ \exists k_{i}\equiv 2k_{j} (mod~ n),~ k_{i}\in S_{nl}^{m-2}  \}   = \\
                      \{    2k_{r} +b _{j}^{m-3}\vert ~  k_{r}\in  S_{nj}^{m-3} \}   $$
            and from it one can again pick out a uniformly distributed set.
            Finally, the remainder will have the view $\{b_{j}^{m'}\}$  where $m'\leq  m-1$,
            thus it will be distributed  uniformly too.      \rule{2mm}{2mm}
\\
\par    {\it Proof of Lemma 1 for the prime  $n>2$.} Let
$\{k _{i}\}=S _{n}$  and  the
            structures of matrices $A$ and $B $ coincide, let $A=(a _{ik} =
            \zeta^{k_{i}}  ),~ B=(b_{ik}  =\zeta^{2k_{i}}   )$. Then $(A ^{n-m}  B^{m} )_{+} =
            E \zeta^{c}\sum_{j=1}^{C_{n}^{m}}\zeta^{a_{j}}  $,
            where $c$ is a constant and $a_{j} =b_{j}^{m}$. According to Theorem  2,  $a_{j}$   are  distributed
            uniformly over the classes of $S_{n}$; then $(A^{n-m}
            B^{m} )_{+} =E  \zeta^{c}\sum_{0}^{n-1} \zeta^{i}/n=0$.
            \rule{2mm}{2mm}
            \\

\par       {\it Proof of Lemma  2.}  Let $ \{A _{i}\vert~ i=1,\ldots,m
\}$ be  concerted, $B_{i} = A _{l}\times A _{i},~ i=1,\ldots,m,~
B_{m+1}   =A_{l}\times E, ~l\not= 1$. We prove, that $\{B_{i}
\vert ~ i=1,\ldots,m+1\}$
            is also  concerted. In
$$                          (B_{1}^{n_{1}}   B _{2}^{n_{2}}\ldots)_{+} ,$$
            including $B_{m+1}$, we divide terms into groups  with  the  same
            right part relative to product sign $\times$; every such a group
            contains $C_{n}^{s}$  different terms with the left part being
           $ P(A_{1}^{n-s}   A_{l}^{s} )$, where $P$ is some permutation. Hence, the sum of
            terms of the group is $(A _{1}^{n-s}  A_{l}^{s} )_{+}\times c=0$.
                        \rule{2mm}{2mm}
                        \\
\par The order of constructed matrices may be reduced due to
the fact that  the set $\{B_{i} =A_{l} \times A'_{i},~ B
_{l}=A_{l}\times E'\vert~  l\not= 1\}$ is concerted if
$\{A_{i}\vert~ i=1,\ldots,k\} $ and $\{A'_{j}\vert ~ j=1,\ldots,
k'\} $ are  concerted too; the orders of $A_{i}$  and $A'_{j}$ may
be different (compare with \cite{8}).
\\
\par            {\it  Proof of Lemma 3 for the prime  $n>2$.}
Let specifically $A=(a _{i,i+1}    = \zeta^{k_{i}} ),~
B=(b_{i,i+1} = \zeta^{2k_{i}}  )$, where $k _{i}\equiv i(mod~ n)$,
the indices are  also  residues  modulo  $n.$  We  prove
            that $(A,B,AB)$ is concerted. One has $AB=\zeta BA$, as

$$ AB=(c_{i,i+2} =  \zeta^{k_{i}+2k_{i+1}}      ),~  BA=(c'_{i,i+2} = \zeta^{2k_{i}+k_{i+1}}        ). $$
            Instead of
$$   (A ^{n_{1}} B^{n_{2}}  (AB)^{n_{3}}  )_{+},~ n_{1} +n_{2} +n_{3} =n,$$
            one can consider the expression
$$  (A ^{n_{1}} B^{n_{2}}  (AB)^{n_{3}}  )_{+} +(A^{n_{1}}  B ^{n_{2}} (BA)^{n_{3}}  )_{+} ,$$
            turning into zero simultaneously with the first one. In the
            last, all terms can be divided into groups, in which $n_{3}$
            numbers of $k_{i}$  from the variable set of $(n_{2} +n_{3} )$
            elements are fixed,
            and, besides, every group has $C_{n}^{n_{2}}$   terms of the view
$$               \zeta^{c}\zeta^{\sum_{1}^{n_{2}} k''_{i}}   ,                  $$
            where $c=\sum_{1}^{n+n_{3}}k'_{i} $ and is not changed, $k'_{i}\in\{ k_{i}\} ,\{ k''_{i}\}= S_{n}$
            and  $ \sum_{1}^{n_{2}} k''_{i} $  are  distributed  uniformly by Theorem 2. Thus
            $(A^{n_{1}}  B^{n_{2}}  (AB)^{}n_{3}  )_{+} =0$.
            \rule{2mm}{2mm}
            \\
\par       {\it Proof of Theorem 1 for the prime $ n>2$.}
The validity of part A of the theorem for the prime $n>2$ is
secured by Lemmas 1  and 2. In proofs of Lemmas it is $A_{i}^{n}
=E $ by construction. We prove now that part B is true. Let $r
_{1}=k$ and the condition (4) be satisfied. \par Rewrite (3) as
$A\Psi=B\Psi$, where $A=\sum_{1}^{k} A_{i} x_{i},~ B=Ey$. \par As
$AB-BA=0$, therefore, one gets $A ^{n}\Psi=B^{n}\Psi $,
multiplying by $A^{n-1}$ from the left, where we have $A ^{n}=E
\sum_{1}^{k}   x_{i}^{n} $  by  (4). \par  As $|A| \not=0$, by the
equalities $\vert A^{n-1}\vert \vert A-B\vert=\vert A^{n} -B
^{n}\vert=0$ it turns out, that if (2) takes place, a determinant
of the system (3)  is  equal to zero. \rule{2mm}{2mm}

\section[3] { The particular case of $n=3$              }

               For $n=3$ the following representation of three concerted
            matrices can be taken:

\begin{displaymath}
  A_{1}=\left(\matrix{ 0 & \zeta & 0 \cr
                                        0 &    0    & \zeta^{2} \cr
                                        1  &     0   &  0  \cr  }
                                        \right),~~
A_{2}=\left(\matrix{   0 &    1    & 0 \cr
                                        0 &    0    & 1 \cr
                                       1  &     0   &  0  \cr  }
                                       \right),~~
A_{3}=\left(\matrix{  0 & 0 & 0 \cr
                                       \zeta  &    0    & 0 \cr
                                       0  & \zeta^{2}  &  0  \cr  }
                                       \right).
\end{displaymath}
Then solutions of the equation $x_{1}^{n} +x_{2}^{n} +x_{3}^{n} =y
^{n}$ can be parametrized with complex quantities $\Psi_{i}$ (
$\delta $ is a determinant of the system of linear equations
relative to $x_{i}$ ) in such the manner:
$$ y\equiv \delta=\Psi_{1}\Psi_{3}^{2}(1-\zeta^{2} )-
\Psi_{3}\Psi_{2}^{2}(1-\zeta)- \Psi_{2}\Psi_{1}^{2} (\zeta-\zeta^{2} ); $$
$$ x_{1} =\zeta\Psi_{1}^{3}  + \zeta^{2}\Psi_{2}^{3} +\Psi_{3}^{3}; $$
$$ x_{2} =-( \zeta\Psi_{1}^{3} +\Psi_{2}^{3} + \zeta^{2}\Psi_{3}^{3} ); $$
$$ x _{3}= \Psi_{3}\Psi_{1}^{2} (1-\zeta^{2} )-
\Psi_{1}\Psi_{2}^{2}(1-\zeta)- \Psi_{2}\Psi_{3}^{2} (\zeta-\zeta^{2} ).$$
\\
\par If $\Psi_{i}$  are integer numbers of the 3-circular field
$K_{3}$, $ x _{i}$ and $y$ should be the same ones. \par This
example shows that search for integer solutions of (1), after
parametrization of it with the help of quantities $\Psi_{i}$,
meets with another problem: it is necessary to study a possibility
of such a choice of complex parameters, that $\Im x_{i} =\Im y_{j}
=0$. \par For (1), Euler's complete rational parametrization for
$n=3 $ is known (see \cite{13}). Also Ramanujan has considered
            this problem \cite{14}. Our parametrization is a general one:
            for any complex values of $x_{i},~ y $ (not only for integer or
            rational), for which the equation is true, a set of
parameters $\Psi_{j}$  exists. But there are similar problems for
both parametrizations: it is difficult to ascertain for which
values of the parameters we will find integer solutions.

\section[4] { Remarks on a solvability of Eqs.~(1).    }
               A combinatorial consideration leads to the restriction:
\begin{theorem}   $$r(n) \leq 2n+1.$$
\end{theorem}
\par          {\it Proof.} Let us consider a set $\{(a_{1},~a
_{2},\ldots,a_{k} )\}$ with $0\leq a_{i}\leq A,~~ a_{i}\leq
a_{i+1}$. The number of elements of the set \cite{15} is

$$                             N(A,k)=C_{A+k}^{k}   =(A+k)!/k!A!,     $$
and the number of different non-zero values of
$\sum_{1}^{k}a_{i}^{n}$   is $kA ^{n}$. If for $k=k_{0}$  some $A$
exists such that $N(A,k_{0})> k_{0} A ^{n}$, then $r(n)\leq 2k_{0}
-1$. As $(A+k)!/A!>A^{k} $, and for $ A^{k_{0}}  >k_{0}!k_{0}
A^{n} $, the inequality $N(A,k _{0})>k_{0} A^{n} $ will be
fulfilled; $k_{0}$  is fixed, therefore $k_{0} =n+1$ is enough.
\rule{2mm}{2mm} \par The bound $r(n)\leq 2n-1$ has been recently
obtained by Habsieger \cite{16}; I am grateful to the referee who
indicated me this fact.
\\
            \par   There are two conjectures about $r(n)$:
$$ 1.~ r(n)\leq n,~ n>1.$$
$$ 2.~ r(n)>e \ln(n),   $$
where $e>0$ is fixed.
\par In my unpublished work \cite{7}, the conjecture:
            $r(p)=p$ for the prime $p$ was mentioned,
            but it is not correct  because of the
            particular solution for $n=5$ \cite{17} (see below).
            \par   I give also the identities which are just for arbitrary
           numbers $a_{i}$  (the proof was given in \cite{7}):
$$ \sum_{k=0}^{n}\sum_{j=1}^{C_{n}^{k}} (g _{jk} ) ^{n}(-1)^{k}\equiv  2^{n} n! \prod_{1}^{n}a_{i},~~
      $$
where $g_{jk}  =\sum_{i=1}^{n}\sigma_{jki}a_{i} $ and
$\sigma_{jki}$ is the element of a sign matrix: $
\vert\sigma_{jki}\vert =1,~~ \sum_{i=1}^{n}\sigma_{jki} =n-2k,~~
\sigma_{jki}= \sigma_{lmi}\delta_{j}^{l}\delta_{k}^{m}$. The
identities have $ 2^{n-1} $  different terms on the left side .
\par   Let us consider now the restrictions on $r(n)$ for $n\leq12$,
obtained by the author under search of particular  solutions of
(1) with  PC AT 386.  I describe briefly the used algorithm which
is not connected  with  linearization of  Eqs. (1). It is based on
the idea of the most compact filling or taking up of some priming
volume $V_{0}$  by single n-cubes, while $V_{0}$  is picked out by
means of  running  over  on  a few cycles. If  $  L \equiv
\sum_{1}^{m} x_{i}^{n},~~ R\equiv \sum_{1}^{k} y_{j}^{n},~~ V =L-R
$ and $ x_{i}\in  [0,A]$ with the fixed $A$, $ y _{j}\in [f
_{1}(L- \sum_{1}^{j-1} y_{l}^{n}  ), f_{2} (L- \sum_{1}^{j-1}
y_{l}^{n}  ) ],$ I shall  refer to such an algorithm as $LmRk$,
for instance $L1R3 $ etc.  Comparing $\vert\vert V_{0}\vert
-y^{n}\vert$   with $ \vert\vert V_{0}\vert -(y+1)^{n}\vert $,  we
choose the smallest of them; then a value  of  $y$  or $ y+1$,
corresponding to it, will be chosen at the first step as a new
element, $ V_{0}$ will  be   replaced   by $\vert V
_{0}\vert-y^{n} $  or $\vert V_{0}\vert -(y+1)^{n} $ and the
procedure will be repeated. If after a fixed number of iterating a
final volume is not equal to zero, $V_{0}$  must be replaced.
Functions $f _{1}$ and $ f_{2}$  were chosen empirically, usually
as $ ((L- \sum_{1}^{j-1} y_{l}^{n}  )/c)^{1/n}  $   with some $c$.
I do not describe the technical problem of representation of big
integers by computations.
            \par   Instead of the evident record
$$  \sum_{1}^{r_{1}} x_{i}^{n}  = \sum_{1}^{r_{2}} y_{j}^{n}$$
let us describe solutions of (1) as
$$ (x_{1} ,x_{2} ,\ldots,x_{r_{1}};y _{1},y_{2} ,\ldots,y_{r_{2}})^{n} =0.            $$
\par It is known that $r(3)=3$. In \cite{7} the solutions were
adduced for $n=4,5,6:~  (3,5,8;7,7)^{4} =0,~~
(4,10,20,28;3,29)^{5} =0,~~ (3,19,22;10,15,23)^{6} =0$. The
referee has indicated me Euler's solution for $n=4:~
(133,134;158,59)^{4} =0$ \cite{13} and the very important one for
$n=5:~ (27,84,110,133;144)^{5} =0 $ \cite{17}. Hence,  it  follows
that $ r(4)\leq 3,~ r(5)\leq 4,~~ r(6)\leq 5$.
\par   I give the best restriction on $r(n)$ for $n=7$:
$$ (149,123,14,10;146,129,90,15)^{7} =0,~ r(7)\leq 7,~~ L3R3;     $$
and the same for$ n=8$:
$$   (43,20,11,10,1;41,35,32,28,5)^{8} =0,~ r(8)\leq 9,~~ L1R3.     $$
For $n=9$ we have:
$$ (73,38,29,9,1;68,67,45,21,18,11,6,4)^{9} =0,~ r(9)\leq 12,~~ L2R1. $$
\par Restrictions on $r(n)$ for $n=10,~11$  are weaker, and the ones are
the best:
$$ (149,42,37,30,25,20,8,5;145,128,100,73,48,13,6,1)^{10}  =0,~
r(10)\leq 15,~~ L2R1;$$
$$(18,6,6,6,4,4,4;17,16,15,13,13,10,9,9,8,1,1,1,1)^{11}  =0,~~
r(11)\leq 19 $$ ($V_{0}$  was chosen by randomization).
\par  Thus, the particular solutions show that $r(n)\leq n \vert~ n\leq7$.

\section[5] { Conclusion                      }

\par  A possibility of  described  linearization  of  Eqs.(1)
sets a number of algebraic problems: the study of algebras of
matrices $\{ A _{i}\}$, construction  of  concerted  sets  of
matrices with the smallest dimension (see \cite{8}) etc.
\par   New quantities $\Psi$ for $n>2$, also as spinors for $n=2$,
permit  to  build  new  algebraic  objects   with   different
transformation laws by transformations of coordinates
$\{x_{i},~y_{j}\} $; research of such objects will be probably of
a great interest.
\par A fundamental fact for the  Diophantine
equations  (1) can be the possibility of parametrization of their
solutions.
\section[6]{ Acknowledgements                    }
I am very grateful to  the  referee  of  my  paper  for
            substantial and useful remarks and suggestions on the
            manuscript which were used by the author to revise
            the manuscript. I must mark with a gratitude that
            the references \cite{9,13,14,16,17} were given me by the referee.

\end{document}